\documentclass[preprint,12pt]{elsarticle}

\usepackage{amssymb}\usepackage{MnSymbol}
\usepackage{amsthm}\usepackage{amsmath}

\journal{Journal of Multivariate Analysis}

\begin{document}

\theoremstyle{plain}
\newtheorem{theorem}{Theorem}
\newtheorem{lemma}[theorem]{Lemma}
\newtheorem{proposition}[theorem]{Proposition}
\newtheorem{corollary}[theorem]{Corollary}
\newtheorem{fact}[theorem]{Fact}
\newtheorem*{main}{Main Theorem}

\theoremstyle{definition}
\newtheorem{definition}[theorem]{Definition}
\newtheorem{example}[theorem]{Example}

\theoremstyle{remark}
\newtheorem*{conjecture}{\indent Conjecture}
\newtheorem{remark}{Remark}
\newtheorem{claim}{Claim}

\newcommand{\N}{\mathbb{N}}
\newcommand{\R}{\mathbb{R}}
\newcommand{\I}{\mathbb{I}}
\newcommand{\B}{\mathbb{B}}
\newcommand{\Vspace}{\vspace{2ex}}
\newcommand{\bfx}{\mathbf{x}}
\newcommand{\bfy}{\mathbf{y}}
\newcommand{\bfh}{\mathbf{h}}
\newcommand{\bfe}{\mathbf{e}}
\newcommand{\Q}{Q}
\newcommand{\os}{\mathrm{os}}
\newcommand{\dd}{\,\mathrm{d}}

\begin{frontmatter}



\title{Extensions of system signatures to dependent lifetimes: Explicit expressions and interpretations}


\author[1]{Jean-Luc Marichal}
\ead{jean-luc.marichal[at]uni.lu}

\author[1]{Pierre Mathonet}
\ead{pierre.mathonet[at]uni.lu}

\address[1]{Mathematics Research Unit, FSTC, University of Luxembourg\\ 6, rue Coudenhove-Kalergi, L-1359 Luxembourg, Grand Duchy of Luxembourg.}

\begin{abstract}
The concept of system signature was introduced by Samaniego for systems whose components have i.i.d.\ lifetimes. We consider its extension to
the continuous dependent case and give an explicit expression for this extension as a difference of weighted means of the structure function
values. We then derive a formula for the computation of the coefficients of these weighted means in the special case of independent continuous
lifetimes. Finally, we interpret this extended concept of signature through a natural least squares approximation problem.
\end{abstract}


\end{frontmatter}



\section{Introduction}

Consider an $n$-component \emph{semicoherent system}. The design of such a system can be described through its \emph{structure function}
$\phi\colon\{0,1\}^n\to\{0,1\}$, which expresses at any time the state of the system in terms of the states of its components. Here
``semicoherent'' means that the structure function $\phi$ is nondecreasing in each variable and satisfies the boundary conditions
$\phi(\mathbf{0})=0$ and $\phi(\mathbf{1})=1$.

Let $X_1,\ldots,X_n$ denote the component lifetimes and let $X_{1:n},\ldots,X_{n:n}$ be the order statistics obtained by rearranging the
variables $X_1,\ldots,X_n$ in ascending order of magnitude; that is, $X_{1:n}\leqslant\cdots\leqslant X_{n:n}$.

The concept of \emph{signature} was introduced in 1985 by Samaniego \cite{{Sam85}} for systems whose components have i.i.d.\ lifetimes as the
$n$-tuple $\mathbf{s}=(s_1,\ldots,s_n)$ defined by
\begin{equation}\label{eq:sa9fd7}
s_k=\Pr(T=X_{k:n})
\end{equation}
where $T$ denotes the system lifetime. That is, $s_k$ is the probability that the $k$th component failure causes the system to fail. For general
background, see Samaniego~\cite{Sam07}.

Under the i.i.d.\ assumption, the probability $s_k$ can be computed as the ratio of $n_k$, the number of orderings for which the $k$th component
failure causes the system failure, to $n!$, the total number of possible orderings of the failure times. An alternative way to calculate $s_k$,
which does not require the exhaustive inspection of all the orderings, was found by Boland \cite{Bol01} through the formula
\begin{equation}\label{eq:asad678}
s_k=\frac{1}{{n\choose n-k+1}}\sum_{\textstyle{\bfx\in\{0,1\}^n\atop |\bfx|=n-k+1}}\phi(\bfx)-\frac{1}{{n\choose
n-k}}\sum_{\textstyle{\bfx\in\{0,1\}^n\atop |\bfx|=n-k}}\phi(\bfx)
\end{equation}
where $|\bfx|=\sum_{i=1}^n x_i$.

Even though the i.i.d.\ assumption is quite reasonable in many applications, especially when we want to compare different system designs, it is
legitimate to investigate the probability (\ref{eq:sa9fd7}) in the general case of dependent lifetimes and to search for formulas which extend
(\ref{eq:asad678}) to this general framework. We observe that only a few results have been obtained in this direction, assuming for instance
that the component lifetimes are exchangeable or independent and exponentially distributed (see, e.g., Samaniego~\cite[{\S}8.3]{Sam07} and
Navarro et al.~\cite{NavSamBalBha08}).

Equation (\ref{eq:asad678}) shows that, in the i.i.d.\ case, the $n$-tuple $\mathbf{s}$ depends only on the system design. This is no longer
true in the general dependent case, where the probability $\Pr(T=X_{k:n})$ may also depend on the joint c.d.f.\ of the component lifetimes, that
is the function
$$
F(t_1,\ldots,t_n)=\Pr(X_1\leqslant t_1,\ldots,X_n\leqslant t_n).
$$
In this general setting, we shall denote the probability $\Pr(T=X_{k:n})$ by $p_k$ to emphasize that the $n$-tuple $\mathbf{p}=(p_1,\ldots,p_n)$
is not always a signature in the strict sense (i.e., a feature of the system design only).


In this paper, assuming only that the joint c.d.f.\ $F$ is absolutely continuous, we provide a closed-form expression for $p_k=\Pr(T=X_{k:n})$ 
as a difference of two weighted arithmetic means of the structure function values whose weights depend only on $F$ (Theorem~\ref{thm:sd98sd7}).
We show that this expression reduces to (\ref{eq:asad678}) as soon as the component lifetimes are exchangeable (Remark~\ref{rem:s9d7f}). We also
provide a useful expression for the weights (as a one-dimensional integral) in the special case of independent lifetimes
(Proposition~\ref{lemma:6d5f67}) and examine the particular case of independent Weibull lifetimes, which includes the exponential model
(Corollary~\ref{cor:vx6cx6c5}). Finally, we show that the $n$-tuple $\mathbf{p}$ can be obtained from a symmetric approximation of the structure
function in the sense of weighted least squares and we point out a formal analogy between this approximation and the concept of projected system
recently introduced in Navarro et al.~\cite{NavSpiBal10}.

Through the usual identification of the elements of $\{0,1\}^n$ with the subsets of $[n]=\{1,\ldots,n\}$, a pseudo-Boolean function
$f\colon\{0,1\}^n\to\R$ can be equivalently described by a set function $v_f\colon 2^{[n]}\to\R$. We simply write $v_f(S)=f(\mathbf{1}_S)$,
where $\mathbf{1}_S$ denotes the $n$-tuple whose $i$th coordinate is $1$, if $i\in S$, and $0$, otherwise. To avoid cumbersome notation, we
henceforth use the same symbol to denote both a given pseudo-Boolean function and its underlying set function, thus writing
$f\colon\{0,1\}^n\to\R$ or $f\colon 2^{[n]}\to\R$ interchangeably.

The $k$th order statistic function $\mathrm{os}_{k:n}\colon\{0,1\}^n\to\{0,1\}$ is defined by $\mathrm{os}_{k:n}(\bfx)=x_{k:n}$. We then have
$\mathrm{os}_{k:n}(\bfx)=1$, if $|\mathbf{x}|\geqslant n-k+1$, and $0$, otherwise. As a matter of convenience, we also formally define
$\mathrm{os}_{0:n}\equiv 0$ and $\mathrm{os}_{n+1:n}\equiv 1$.

\section{Explicit expressions}

Let $F$ be the (absolutely continuous) joint c.d.f.\ of the component lifetimes $X_1,\ldots,X_n$. We define the associated \emph{relative
quality function} $q\colon 2^{[n]}\to [0,1]$ as
$$
q(S) = \Pr\Big(\max_{i\in [n]\setminus S}X_i < \min_{j\in S}X_j\Big)
$$
with the convention that $q(\varnothing)=q([n])=1$. That is, $q(S)$ is the probability that the lifetime of every component in $S$ is greater
than the lifetime of every component in $[n]\setminus S$. Thus defined, $q(S)$ is a measure of the overall quality of the components in $S$ when
compared with the components in $[n]\setminus S$.

Since the r.v.'s $X_1,\ldots,X_n$ are continuous, we see that the function $q$ can also be written as
\begin{equation}\label{eq:sadf65sdf}
q(S)=\sum_{\sigma\in\mathfrak{S}_n\, :\, \{\sigma(n-|S|+1),\ldots,\sigma(n)\}=S}\Pr(X_{\sigma(1)}<\cdots <X_{\sigma(n)})
\end{equation}
where $\mathfrak{S}_n$ denotes the group of permutations on $[n]$.

We then observe that, for every $k\in [n]$, the values $q(S)$ for $|S|=k$ sum up to one. In fact, by (\ref{eq:sadf65sdf}), we have
\begin{equation}\label{eq:9s8d}
\sum_{\textstyle{S\subseteq [n]\atop |S|=k}}q(S) = \sum_{\sigma\in \mathfrak{S}_n}\Pr(X_{\sigma(1)}<\cdots <X_{\sigma(n)}) = 1.
\end{equation}

\begin{remark}\label{rem:a9sdf76}
\begin{enumerate}
\item[(a)] The validity of (\ref{eq:9s8d}) is especially transparent if one focuses on a particular simple case, say, when $k = 1$. It then
simply says that there exists (with probability $1$) a unique component $j$ whose lifetime is maximum, that is,
$$
\sum_{j=1}^n \Pr\Big(\max_{i\in [n]\setminus\{j\}}X_i<X_j\Big)=1.
$$

\item[(b)] If the variables $X_1,\ldots,X_n$ are exchangeable, then the relative quality function $q$ is symmetric. In this case, by
(\ref{eq:9s8d}), we have $q(S)=1/{n\choose |S|}$.

\item[(c)] Equation~(\ref{eq:9s8d}) shows that comparing $q(S)$ with $q(S')$ is relevant whenever $|S|=|S'|$. In general, according to (b)
above, the relative quality could be better measured by the function $\tilde{q}(S)={n\choose |S|}\, q(S)$.
\end{enumerate}
\end{remark}

We now give an expression for $p_k=\Pr(T=X_{k:n})$ as a difference of two arithmetic means of the structure function values weighted by the
relative quality function. We first present a lemma.

\begin{lemma}\label{lemma:7sf987}
For every $k\in [n]$, we have
$$
\Pr(T\geqslant X_{k:n})=\sum_{|\bfx|=n-k+1}q(\bfx)\,\phi(\bfx).
$$
\end{lemma}

\begin{proof}
For every $k\in [n]$ and every $\sigma\in \mathfrak{S}_n$, we have
\begin{equation}\label{eq:sf98sdf}
\Pr(T\geqslant X_{\sigma(k)}\mid X_{\sigma(1)}<\cdots <X_{\sigma(n)})=\phi\big(\{\sigma(k),\ldots,\sigma(n)\}\big).
\end{equation}
Indeed, assume $X_{\sigma(1)}<\cdots <X_{\sigma(n)}$. Using the path representation of the life function \cite{EsaMar70} of the system, we have
$T\geqslant X_{\sigma(k)}$ if and only if
\begin{equation}\label{eq:9sdf8d}
\max_{1\leqslant j\leqslant m}\;\min_{i\in P_j} X_i\geqslant X_{\sigma(k)}
\end{equation}
where $P_1,\ldots,P_m$ denote the minimal path sets of the system. In turn, event (\ref{eq:9sdf8d}) is equivalent to
$$
\mbox{``There is a minimal path set $P$ such that $X_i\geqslant X_{\sigma(k)}$ for all $i\in P$''}
$$
Equivalently, ``there is a minimal path set $P\subseteq \{\sigma(k),\ldots,\sigma(n)\}$''. By monotonicity of $\phi$, this event reduces to
``$\phi\big(\{\sigma(k),\ldots,\sigma(n)\}\big)=1$'', which finally proves (\ref{eq:sf98sdf}).

By combining the law of total probability with (\ref{eq:sf98sdf}), we obtain
$$
\Pr(T\geqslant X_{k:n}) = \sum_{\sigma\in\mathfrak{S}_n}\phi(\{\sigma(k),\ldots,\sigma(n)\})\,\Pr(X_{\sigma(1)}<\cdots <X_{\sigma(n)}).
$$
Grouping the terms for which $\{\sigma(k),\ldots,\sigma(n)\}$ is a fixed set $S$ and then summing over $S$, we obtain
$$
\Pr(T\geqslant X_{k:n}) = \sum_{|S|=n-k+1}~\sum_{\sigma\in\mathfrak{S}_n\, :\,
\{\sigma(k),\ldots,\sigma(n)\}=S}\phi(S)\,\Pr(X_{\sigma(1)}<\cdots <X_{\sigma(n)}).
$$
The result then follows from (\ref{eq:sadf65sdf}).
\end{proof}

\begin{theorem}\label{thm:sd98sd7}
For every $k\in [n]$, we have
\begin{equation}\label{eq:sdf67}
p_k=\sum_{|\bfx|=n-k+1}q(\bfx)\,\phi(\bfx)-\sum_{|\bfx|=n-k}q(\bfx)\,\phi(\bfx).
\end{equation}
\end{theorem}

\begin{proof}
We have $\Pr(T=X_{k:n})=\Pr(T\geqslant X_{k:n})-\Pr(T\geqslant X_{k+1:n})$. We then conclude by Lemma~\ref{lemma:7sf987}.
\end{proof}

\begin{remark}\label{rem:s9d7f}
\begin{enumerate}
\item[(a)] It is noteworthy that $p_k$ can be rewritten in the form
$$
p_k=\sum_{\bfx\in\{0,1\}^n}r_k(\bfx)\,\phi(\bfx),
$$
with $r_k(\bfx)=q(\bfx)\,(-x_{k+1:n}+2\, x_{k:n}-x_{k-1:n})$. This fact follows immediately from the identity
$$
\sum_{\bfx\in\{0,1\}^n} q(\bfx)\,\phi(\bfx)\, x_{k:n}=\sum_{|\bfx|\geqslant n-k+1} q(\bfx)\,\phi(\bfx).
$$


\item[(b)] As expected, we observe from (\ref{eq:sdf67}) that $p_k$ combines linearly two pieces of information:
\begin{enumerate}
\item[(i)] the system design, which is encoded in the structure function $\phi$, and

\item[(ii)] the component lifetimes, which are encoded in the relative quality function $q$.
\end{enumerate}

\item[(c)] When the variables $X_1,\ldots,X_n$ are exchangeable, by Remark~\ref{rem:a9sdf76}$(b)$ we see that (\ref{eq:sdf67}) reduces to
(\ref{eq:asad678}). Thus in this case, $\mathbf{p}$ is the signature of the system. This fact was previously observed in \cite[Lemma
1]{NavRyc07}.
\end{enumerate}
\end{remark}

The following proposition gives a formula for $q(S)$ as a one-dimensional integral in the special case of independent lifetimes.

\begin{proposition}\label{lemma:6d5f67}
For independent continuous lifetimes $X_1,\ldots,X_n$, each $X_i$ having p.d.f.\ $f_i$ and c.d.f.\ $F_i$, with $F_i(0)=0$, we have
$$
q(S)=\sum_{j\in S}\int_0^{\infty} f_j(t)\,\prod_{i\in [n]\setminus S}F_i(t)\,\prod_{i\in S\setminus\{j\}}(1-F_i(t))\dd t.
$$
\end{proposition}

\begin{proof}
Denote the p.d.f.\ and c.d.f.\ of a r.v.\ $X$ by $f_X$ and $F_X$, respectively. Recall that for two independent continuous r.v.'s $X$ and $Y$,
we have
$$
\Pr(X<Y) = \iint_{x<y}f_X(x)f_Y(y)\dd x\dd y = \int_{-\infty}^{\infty}F_X(y)f_Y(y)\dd y.
$$
By applying this formula to $q(S)$, we obtain
\begin{eqnarray*}
q(S) &=& \int_0^\infty\Pr\Big(\max_{i\in [n]\setminus S}X_i \leqslant t\Big)\,\frac{\mathrm{d}}{\mathrm{d}t}\Pr\Big(\min_{j\in S}X_j\leqslant t\Big)\dd t\\
&=& \int_0^\infty\prod_{i\in [n]\setminus S}F_i(t)\,\frac{\mathrm{d}}{\mathrm{d}t}\bigg(1-\prod_{i\in S}(1-F_i(t))\bigg)\dd t
\end{eqnarray*}
which immediately leads to the result.
\end{proof}

\begin{corollary}\label{cor:vx6cx6c5}
For independent Weibull lifetimes, with $F_i(t)=1-e^{-(\lambda_i t)^{\alpha}}$, 
we have
\begin{equation}\label{eq:ycx7v6}
q(S)=\sum_{K\subseteq [n]\setminus S}(-1)^{|K|}\,\frac{\lambda_{\alpha}(S)}{\lambda_{\alpha}(K\cup S)}
\end{equation}
for every $S\neq\varnothing$, where $\lambda_{\alpha}(S)=\sum_{i\in S}\lambda_i^{\alpha}$.
\end{corollary}

\begin{proof}
By Proposition~\ref{lemma:6d5f67}, we have
$$
q(S)=\sum_{j\in S}\lambda_j^{\alpha}\,\int_0^{\infty}{\alpha}\, t^{{\alpha}-1}\,e^{-\lambda_{\alpha}(S)\, t^{\alpha}}\,\prod_{i\in [n]\setminus
S} \big(1-e^{-\lambda_{\alpha}(\{i\})\, t^{\alpha}}\big)\dd t
$$
where the product can be expanded (by the generalized binomial theorem) as
$$
\prod_{i\in [n]\setminus S} \big(1-e^{-\lambda_{\alpha}(\{i\})\, t^{\alpha}}\big)=\sum_{K\subseteq [n]\setminus S} (-1)^{|K|}\,
e^{-\lambda_{\alpha}(K)\, t^{\alpha}}.
$$
We then have
$$
q(S)=\sum_{j\in S}\lambda_j^{\alpha}\sum_{K\subseteq [n]\setminus S} (-1)^{|K|}\int_0^{\infty}{\alpha}\, t^{{\alpha}-1}\,
e^{-\lambda_{\alpha}(K\cup S)\, t^{\alpha}}\dd t
$$
which immediately leads to the result.
\end{proof}

\begin{remark}
\begin{enumerate}
\item[(a)] Given a set function $q\colon 2^{[n]}\to [0,1]$, there is an additive set function $\lambda_{\alpha}\colon
2^{[n]}\to\left]0,\infty\right[$ satisfying (\ref{eq:ycx7v6}) if and only if
\begin{equation}\label{eq:asd9f7x}
q([n]\setminus\{i\})>0\qquad \forall i\in [n]
\end{equation}
and
\begin{equation}\label{eq:asd9f7}
q(S) = \sum_{K\subseteq [n]\setminus S}(-1)^{|K|}\,\frac{\sum_{i\in S}q([n]\setminus\{i\})}{\sum_{i\in K\cup S}q([n]\setminus\{i\})}\qquad
\forall S\neq\varnothing.
\end{equation}
Indeed, if such an additive function exists, then by (\ref{eq:ycx7v6}) we obtain
\begin{equation}\label{eq:asd9f7a}
q([n]\setminus\{i\})=1-\frac{\lambda_{\alpha}([n]\setminus\{i\})}{\lambda_{\alpha}([n])}=\frac{\lambda_{\alpha}(\{i\})}{\lambda_{\alpha}([n])}
\end{equation}
which leads to (\ref{eq:asd9f7x}) and (\ref{eq:asd9f7}). Conversely, if (\ref{eq:asd9f7x}) holds, then we may choose
$\lambda_{\alpha}(\{i\})=q([n]\setminus\{i\})$ and we then see that (\ref{eq:asd9f7}) leads to (\ref{eq:ycx7v6}). Thus (\ref{eq:asd9f7x}) and
(\ref{eq:asd9f7}) provide necessary and sufficient conditions on a set function $q\colon 2^{[n]}\to [0,1]$ to be a relative quality function
obtained from Weibull lifetimes (with a common shape parameter $\alpha$).

\item[(b)] Under the assumptions of Corollary~\ref{cor:vx6cx6c5}, by (\ref{eq:asd9f7a}) the ratio
$\lambda_{\alpha}(\{i\})/\lambda_{\alpha}([n])$ is exactly the probability that $X_i$ is the shortest lifetime. More generally, the ratio
$\lambda_{\alpha}(S)/\lambda_{\alpha}([n])=\sum_{i\in S}q([n]\setminus\{i\})$ is the probability that the component having the shortest lifetime
is in $S$.

\item[(c)] For $S\subseteq [n]$, the $S$-difference of a function $f\colon\{0,1\}^n\to\R$ is defined inductively by $\Delta^{\varnothing}f=f$
and $\Delta^Sf=\Delta^{\{i\}}\Delta^{S\setminus\{i\}}f$ for $i\in S$, with $\Delta^{\{i\}}f(\mathbf{x})=f(\mathbf{x}\mid x_i=1)-f(\mathbf{x}\mid
x_i=0)$. It is then easy to see \cite[{\S}2]{GraMarRou00} that (\ref{eq:ycx7v6}) can be rewritten as
\begin{equation}\label{eq:sd6f7}
q(S)=(-1)^{n-|S|}\,\lambda_{\alpha}(S)\,\Big(\Delta^{[n]\setminus S}\frac{1}{\lambda_{\alpha}}\Big)(S).
\end{equation}
Moreover, the $([n]\setminus S)$-difference in (\ref{eq:sd6f7}) can be interpreted as the marginal interaction \cite[{\S}2]{GraMarRou00}
(associated with the function $1/\lambda_{\alpha}$) among the components in $[n]\setminus S$ conditioned to the presence of the components in
$S$.
\end{enumerate}
\end{remark}

\section{Links with approximations of structure functions}

In \cite{MarMat} the authors solved the problem of approximating a given pseudo-Boolean function $f\colon\{0,1\}^n\to\R$ by a symmetric one in
the sense of weighted least squares.

Specifically, given a weight function $w\colon\{0,1\}^n\to\left]0,\infty\right[$, the \emph{best symmetric approximation} of a function
$f\colon\{0,1\}^n\to\R$ is defined as the unique symmetric function $f^*\colon\{0,1\}^n\to\R$ that minimizes the weighted squared distance
\begin{equation}\label{eq:xc66xxy}
\|f-g\|^2=\sum_{\bfx\in\{0,1\}^n}w(\bfx)\big(f(\bfx)-g(\bfx)\big)^2
\end{equation}
among all symmetric functions $g\colon\{0,1\}^n\to\R$.

The best symmetric approximation $f^*$ is actually the orthogonal projection of $f$, with respect to the inner product
$$
\langle f,g\rangle=\sum_{\bfx\in\{0,1\}^n}w(\bfx)f(\bfx)g(\bfx)
$$
onto the linear subspace of symmetric functions $g\colon\{0,1\}^n\to\R$. In terms of the order statistic functions, this projection is given by
\begin{equation}\label{eq:sd987x}
f^*=f(\mathbf{0})+\sum_{k=1}^nc_k\,\mathrm{os}_{k:n}
\end{equation}
where
\begin{equation}\label{eq:sd987}
c_k=\sum_{|\bfx|=n-k+1}\overline{w}(\bfx)\,f(\bfx)-\sum_{|\bfx|=n-k}\overline{w}(\bfx)\,f(\bfx)
\end{equation}
and
$$
\overline{w}(\bfx)=\frac{w(\bfx)}{\sum_{|\mathbf{z}|=|\mathbf{x}|}w(\mathbf{z})}
$$
which shows that $c_k$ is actually a difference of two expected values (see \cite{MarMat}). Since $c_k$ is the coefficient of
$\mathrm{os}_{k:n}$ in $f^*$, it can be interpreted as a measure of the influence of the $k$th smallest variable on $f$.

Now, consider an $n$-component semicoherent system defined by a structure function $\phi$ and an absolutely continuous joint c.d.f.\ $F$ of the
component lifetimes. Assume that the associated relative quality function $q$ is strictly positive. Consider also the weighted distance
(\ref{eq:xc66xxy}) with $w=q$ and apply the approximation problem above to the structure function $f=\phi$. By (\ref{eq:9s8d}) we see that
$\overline{w}=\overline{q}=q$. Theorem~\ref{thm:sd98sd7} then shows that the coefficient $c_k$, as defined in (\ref{eq:sd987}), is precisely the
probability $p_k=\Pr(T=X_{k:n})$.

Moreover, from (\ref{eq:sd987x}) it follows that the best symmetric approximation $\phi^*$ of $\phi$ (with respect to the weighted distance
(\ref{eq:xc66xxy}) with $w=q$) is given by
\begin{equation}\label{eq:sd987y}
\phi^*=\sum_{k=1}^np_k\,\mathrm{os}_{k:n}.
\end{equation}
Under the i.i.d.\ assumption, (\ref{eq:sd987y}) reduces to $ \phi^*=\sum_{k=1}^ns_k\,\mathrm{os}_{k:n}, $ where $s_k$ is given by
(\ref{eq:asad678}).


Interestingly enough, we also observe a formal analogy between the orthogonal projection (\ref{eq:sd987y}) of $\phi$ and the concept of
\emph{projected system}, recently introduced in Navarro et al.~\cite{NavSpiBal10}. Indeed, $\phi^*$ is a combination of the order statistic
functions weighted by the $n$-tuple $\mathbf{p}$ while the projected system is the system which mixes the $k$-out-of-$n$ systems
($\phi_k=\mathrm{os}_{k:n}$) with mixing distribution $\mathbf{p}$.

\section*{Acknowledgments}

The authors wish to thank the reviewers for helpful comments and suggestions. This research is supported by the internal research project
F1R-MTH-PUL-09MRDO of the University of Luxembourg.





\bibliographystyle{elsarticle-num}



\end{document}